\documentclass[reqno]{amsart}

\usepackage{amssymb,amsthm,amsmath,dsfont}

\newtheorem{defi}{Definition}[section]
\newtheorem{thm}[defi]{Theorem}

\newtheorem{cor}[defi]{Corollary}

\newtheoremstyle{normal}% name
    {}              % Space above
    {}              % Space below
    {\normalfont}       % Body font
    {}                  % Indent amount (empty = no indent, \parindent = para indent)
    {\bfseries}         % Thm head font
    {.}                 % Punctuation after thm head
    {.5em}              % Space after thm head: " " = normal interword space; \newline = linebreak
    {}                  % Thm head spec (can be left empty, meaning `normal')
\theoremstyle{normal}

\newtheorem{rem}[defi]{Remark}

\newcommand{\bpr}{\begin{proof}[Proof]}  %Proof
\newcommand{\epr}{\end{proof}}
\newcommand{\beq}{\begin{equation}}
\newcommand{\eeq}{\end{equation}}
\newcommand{\bce}{\begin{center}}
\newcommand{\ece}{\end{center}}
\newcommand{\be}{\begin{enumerate}}  %enumerate
\newcommand{\ee}{\end{enumerate}}

\newcommand{\compemb}{\lhook\hspace{-0.151cm}{-}\hspace{-0.3cm}\hookrightarrow}

\DeclareMathOperator*{\dist}{dist}
\DeclareMathOperator*{\diver}{div}

\def\pa{\partial}
\def\om{\omega}

\def\ep{\varepsilon}
\def\de{\delta}

\def\ga{\gamma}

\def\Om{\Omega}

\def\R{\mathbb R}

\def\N{\mathbb N}

\def\B{\mathbb B}
\def\E{\mathbb E}

\def\calB{\mathcal B}

\def\calT{\mathcal T}

\def\calR{\mathcal R}

\def\calM{\mathcal M}

\def\calU{\mathcal U}

\def\calV{\mathcal V}

\numberwithin{equation}{section}

\def\nn{\nonumber}

\title[Quasilinear parabolic evolution equations in weighted $L_p$-spaces]
{On quasilinear parabolic evolution equations in weighted $L_p$-spaces II}

\author[Jeremy LeCrone]{Jeremy LeCrone}
\address{Department of Mathematics, Kansas State University,
138 Cardwell Hall, Manhattan, KS 66506, USA}
\email{lecronjs@ksu.edu}

\author[Jan Pr\"uss]{Jan Pr\"uss}
\address{Institut f\"ur Mathematik, Martin-Luther-Universit\"at Halle-Wittenberg,
Theodor-Lieser-Str. 5, 06120 Halle, Germany}
\email{jan.pruess@mathematik.uni-halle.de}

\author[Mathias Wilke]{Mathias Wilke}
\address{Institut f\"ur Mathematik, Martin-Luther-Universit\"at Halle-Wittenberg,
Theodor-Lieser-Str. 5, 06120 Halle, Germany}
\email{mathias.wilke@mathematik.uni-halle.de (corresponding author)}


\date{\today}


\begin{document}

\thanks{The work of Jeremy LeCrone was partially supported by a grant from the Simons Foundation (\#245959 to
G. Simonett)}

\maketitle

\begin{abstract}
Our study of abstract quasi-linear parabolic problems in time-weighted $L_p$-spaces, begun in \cite{KPW10}, is extended in this paper to include singular lower order terms, while keeping low initial regularity. The results are applied to reaction-diffusion problems, including Maxwell-Stefan diffusion, and to geometric evolution equations like the surface-diffusion flow or the Willmore flow.
The method presented here will be applicable to other parabolic systems, including free boundary problems.
\end{abstract}

\section{Introduction}

\noindent
In the recent decades, the $L_p$-theory for abstract parabolic evolution equations has become very useful in applications to
parabolic partial differential equations and systems, in particular for parabolic free boundary problems like Stefan problems with surface tension
or two-phase Navier-Stokes problems, and also in the theory of geometric evolution equations, for example, the mean curvature flow, the surface-diffusion flow,
the Willmore flow, and many others. The theory of abstract parabolic evolution equations offers a framework to answer questions concerning local well-posedness,
construction of the corresponding local semi-flow and its qualitative behaviour, like stability of equilibria, invariant manifolds near equilibria, and asymptotic
behaviour of the solutions, especially their convergence to an equilibrium.

There are numerous contributions to this field \cite{Ama88,Ama89,Ama90,Ama93,Ama05,CleSi01,daPraGris79,Esch94,LSU,Lun95,PSZ,Sim94,Sim95,Yag91}; this is just a selection.

Here we follow the approach based on $L_p$-maximal regularity. This approach was first introduced in the paper by Cl\'ement and Li \cite{CleLi93} and has been further developed in \cite{KPW10,MeySchn12a,JanBari}. Following the ideas in \cite{PrSi04}, this approach has recently been extended to weighted $L_p$-spaces in \cite{KPW10} for several reasons.
Time-weights can be used to lower initial regularity and to exploit parabolic regularization which is typical for quasilinear parabolic problems.
This has consequences for compactness properties of solutions, and on the quality of a priori estimates needed for global existence.
As a byproduct, time-weights can also be used to reduce compatibility conditions to a minimum.
The idea of employing time-weights
is by no means new, at has been employed before in papers by \cite{Ang90,Lun95,MeySchn12a} and others.

In this paper we continue our work begun in \cite{KPW10}, extending our results to  equations which have singular lower order terms, while keeping low initial regularity.
To make this statement more precise we begin with the setting of \cite{KPW10}. Let $X_0,X_1$ be Banach spaces such that $X_1\hookrightarrow X_0$ densely, and consider the abstract evolution equation
\begin{equation}\label{1}
\dot{u}+A(u)u=F(u),\; t>0,\quad u(0)=u_0.
\end{equation}
In the framework of $L_p$-theory one looks for solutions $u$ in the class
$$u\in \mathbb{E}_1(J):=H^1_p(J;X_0)\cap L_p(J;X_1)\hookrightarrow C(J;X_\gamma),$$
where $J=[0,T]$ is a finite interval, and  $X_\gamma = (X_0,X_1)_{1-1/p,p}$ is the time trace space of this class. $X_\gamma$ is the natural state space for the semi-flow
generated by \eqref{1}, in the $L_p$-setting, provided it exists. In particular, the initial value $u_0$ must belong to this space, hence one does not see a {\em parabolic regularization}. Choosing $p$ large, the space $X_\gamma$ typically will be a space of {\em high regularity}, and so to obtain a global solution the a priori bounds needed for global existence must necessarily be bounds in $X_\gamma$, i.e.\ they must be of high quality.

To lower the initial regularity - and at the same time the quality of necessary a priori bounds - we employ time-weighted $L_p$-spaces. Defining
\begin{align*}
L_{p,\mu}(J;X_1)&=\{u\in L_{1,loc}(J;X_1):\, t^{1-\mu}u\in L_p(J;X_1)\},\\
 H^1_{p,\mu}(J;X_0)&=\{u\in L_{p,\mu}(J;X_0)\cap H^1_{1}(J;X_0):\, t^{1-\mu}\dot{u}\in L_p(J;X_0)\},
 \end{align*}
we now look for solutions in the class
$$u\in \mathbb{E}_{1,\mu}(J):=H^1_{p,\mu}(J;X_0)\cap L_{p,\mu}(J;X_1)\hookrightarrow C(J;X_{\gamma,\mu}).$$
The trace space for this class is $X_{\gamma,\mu}=(X_0,X_1)_{\mu-1/p,p},$ with $\mu\in (1/p,1]$; note that $X_\gamma=X_{\gamma,1}$, and that
$$ \mathbb{E}_{1,\mu}([0,T])\hookrightarrow \mathbb{E}_1([\delta,T]),$$
for any small $\delta\in (0,T)$, showing instant smoothing of solutions. This striking property of quasilinear parabolic problems allows to pass from bounds on solutions in $X_{\gamma,\mu}$ to bounds in $X_\gamma$, and even to compactness in $X_\gamma$ if the scale $(X_0,X_1)_{\theta,p}$ is a compactly embedded scale.

In our previous paper \cite{KPW10} we worked out this approach under the assumptions
$$(A,F)\in C^{1-}(V_\mu; \mathcal{B}(X_1,X_0)\times X_0),\quad  V_\mu\subset X_{\gamma,\mu} \; \text{open},$$
and $A(u)$ has the property of $L_p$-maximal regularity, for each $u\in V_\mu$. The assumption on $A$ cannot be relaxed substantially, but that on $F$ is not optimal and can be improved.
To explain this, consider the following three scalar parabolic equations in a bounded smooth domain $\Omega\subset\R^n$.
\begin{align}
\partial_t u -a(u) \Delta u &=f(u)\label{examples1}\\
\partial_t u -a(u,\nabla u)\Delta u&= f(u,\nabla u)\label{examples2}\\
\partial_t u - \operatorname{div}(a(u)\nabla u) &=f(u)\label{examples3}.
\end{align}
We equip these problems with the Neumann boundary condition $\partial_\nu u=0$ on $\partial\Omega\in C^2$. To ensure parabolicity, we suppose $a\ge a_0>0$. For all these problems we choose as a base space $X_0=L_q(\Omega)$, where $q\in (1,\infty)$. Then $X_1=\{u\in H^2_q(\Omega):\, \partial_\nu u=0\}$,
$$X_{\gamma,\mu}=\{u\in B_{qp}^{2\mu-2/p}(\Omega):\, \partial_\nu u=0\},$$
whenever the trace $\partial_\nu u$ exists and $X_{\gamma,\mu}=B_{qp}^{2\mu-2/p}(\Omega)$ otherwise; the latter is the case if $2\mu< 1+ 1/q+2/p$. Here $B_{qp}^s$ denote the Besov spaces, see e.g.\ Triebel \cite{Tri83}.

To solve \eqref{examples1} we need $u\in C(J\times\bar{\Omega})$, which follows if $X_{\gamma,\mu}\hookrightarrow C(\bar{\Omega})$, i.e.\ $2\mu > 2/p + n/q$. Such a choice of $\mu \in (1/p,1]$ is possible if $2/p+n/q<2$. Note that no compatibility condition will show up if $q > n-1$ and $\mu$ is chosen sufficiently small!
Next consider \eqref{examples2}; here we need $u\in C(J;C^1(\overline{\Omega}))$ to be able to solve the problem. This means
$ 2\mu > 1 +2/p+n/q$, and so the condition $2/p+n/q<1$ is necessary. In this case there is no way to avoid the trace of $\partial_\nu u$, which means that the compatibility condition
$\partial_\nu u_0=0$ on $\partial\Omega$ must be satisfied.

Finally,  we  look at \eqref{examples3}, which can be rewritten in explicit form as
\begin{equation}\label{challenge}
 \partial_t u -a(u)\Delta u = f(u) +a^\prime(u)|\nabla u|^2.
\end{equation}
Ignoring the second term on the right hand side, the optimal choice for $\mu$ would be as for \eqref{examples1}, i.e.\ $2\mu> 2/p+n/q$. But then, what about the second term on the right hand side, which is quadratic in $\nabla u$? The key observation is that this term is still of lower order {\em and } of polynomial growth (here second order) in $\nabla u$! Such terms ask for estimates of Gagliardo-Nirenberg type (see e.g.\ \cite{Fried69}), and precisely this is what we want to economize in this paper.

To implement this idea we split $F=F_r+F_s$, where $F_r\in C^{1-}(V_\mu;X_0)$ as in our previous paper \cite{KPW10} (subscript $r$ means regular). The singular part
$ F_s:V_\mu\cap X_\beta \to X_0$ with $X_\beta:=(X_0,X_1)_{\beta,p}$ in the simplest case satisfies the following condition:

\medskip

\noindent
{\em There exists numbers $\beta\in (\mu-1/p,1)$, $\rho\geq0$ with
\begin{equation}\label{restriction}
 (1+\rho)(\beta-(\mu-1/p))<1-(\mu-1/p),
\end{equation}
and a continuous function $c:V_\mu\times V_\mu \to \R$ such that for all $u,\bar{u}\in V_\mu\cap X_\beta$
\begin{equation}\label{Fs}
|F_s(u)-F_s(\bar{u})|_{X_0}\leq c(|u|_{X_{\gamma,\mu}},|\bar{u}|_{X_{\gamma,\mu}})\left[
1+|u|_{X_{\beta}}^\rho+|\bar{u}|_{X_{\beta}}^\rho\right]|u-\bar{u}|_{X_\beta}.
\end{equation}}
Here $1+\rho$ is the order of the nonlinearity $F_s$ given by the equation in question. $\beta$ is also given by the equation under consideration, but it depends also on the scale generated by the spaces $X_0$ and $X_1$; in example \eqref{challenge} we will have $1+\rho=2$ and $\beta= 1/2 + n/4q+\varepsilon$ with a small $\varepsilon>0$, see Section \ref{sec:appl} for details. Then \eqref{restriction} defines a lower bound for $\mu>1/p$; for \eqref{challenge} this will be $2\mu>2/p+n/q$, which is already needed for $X_{\gamma,\mu}\hookrightarrow C(\bar{\Omega})$, so in this example there are no further restrictions on $\mu$.

We show that with assumptions \eqref{restriction}, \eqref{Fs}  on the singular part $F_s$ (or its refined version \eqref{LWP:AssF_20}, \eqref{LWP:AssF_2} in Section \ref{LWPsec}) the $L_p$-theory with time-weights for \eqref{1} developed in \cite{KPW10} remains valid. These conditions are flexible enough to cover parabolic systems including \eqref{challenge} as well as many other parabolic problems. In particular, in Section \ref{sec:appl} we show that for the surface diffusion flow and the Willmore flow
it is enough to assume the initial regularity $B^{4\mu-4/p}_{qp}$ with $4\mu> 1+4/p+n/q$. In the setting of the H\"{o}lder spaces $C^\alpha$ this has been recently observed in \cite{Asai12}.

This paper is structured as follows. Section \ref{LWPsec} is concerned with the local well-posedness and optimal regularity of \eqref{1}. To be precise, we will prove the existence and uniqueness of local-in-time solutions in the class $\mathbb{E}_{1,\mu}(0,T)$ depending Lipschitz continuously on the initial data. Furthermore we show that each local solution can be continued to a maximal solution being defined on a maximal interval of existence $J(u_0):=[0,t^+(u_0))$. In Section \ref{sec:relcomp} we will show that if $X_{\gamma,\bar{\mu}}$ is compactly embedded into $X_{\gamma,\mu}$ and if the solution $u(t)$ of \eqref{1} satisfies
$$u\in BC([\tau,t^+(u_0));X_{\gamma,\bar{\mu}})$$
for some $\tau\in (0,t^+(u_0))$ and $\bar{\mu}\in (\mu,1]$, then the set $\{u(t)\}_{t\in[\delta,t^+(u_0))}$ is already relatively compact in $X_{\gamma,1}$ for any $\delta\in (0,t^+(u_0))$. This in turn yields the global existence of the solution.
%Section \ref{sec:regprop} is devoted to the study of the regularity properties of the solution $u(t)$ of \eqref{Int1}. Under the additional assumption that $A$ and $F_j$, $j\in\{1,2\}$ are $k$-times continuously differentiable, where $k\in\mathbb{N}$, we show that the solution is $C^k$ in $t$ and $u_0$.
Finally, in Section \ref{sec:appl} we show how to apply our results to concrete problems. One of these problems is considered with \emph{reaction diffusion systems} in the context of \emph{Maxwell-Stefan diffusion}, which is in the focus of current mathematical research. The other application deals with geometric evolution laws, namely the \emph{surface diffusion flow} and the \emph{Willmore flow}. We consider the situation where each member of a time dependent family of hypersurfaces is given as a graph of a height function over a fixed and bounded domain. This allows one to transform the evolution law for the hypersurfaces to a fourth order quasilinear parabolic equation for the height function to which we can apply the results in this article. In particular we establish an $L_p$--$L_q$-theory.

\medskip

\noindent
\textbf{Notations.} For the spaces
$$\E_{0,\mu}(0,T):=L_{p,\mu}(0,T;X_0),$$
and $\E_{1,\mu}(0,T)$ we use the norms
$$\|u\|_{\mathbb{E}_{0,\mu}(0,T)}=\|u\|_{L_{p,\mu}(0,T;X_0)}:=\|[t\mapsto t^{1-\mu}u(t)]\|_{L_p(0,T;X_0)},$$
and
$$\|u\|_{\mathbb{E}_{1,\mu}(0,T)}:=\|u\|_{\mathbb{E}_{0,\mu}(0,T)}+
\|\dot{u}\|_{\mathbb{E}_{0,\mu}(0,T)}+\|u\|_{L_{p,\mu}(0,T;X_1)},$$
which turn them into Banach spaces. Furthermore we denote by $||\cdot||_{\infty,X_{\gamma,\mu}}$ the norm in $C([0,T];X_{\gamma,\mu})$. In the \emph{classical} case $\mu=1$ we simply use the notation $\E_0$, $\E_1$ and $X_\gamma$ instead of $\E_{1,1}$, $\E_{0,1}$ and $X_{\gamma,1}$. We will also make use of the notation
$$_0\E_{1,\mu}(0,T):=\{u\in\E_{1,\mu}(0,T):u|_{t=0}=0\}.$$
If $M_1$ and $M_2$ are metric spaces and $G:M_1\to M_2$, then $G\in C^{1-}(M_1;M_2)$ means that $G$ is locally Lipschitz continuous. Furthermore we write $X_1{\compemb} X_0$ if $X_1$ is compactly embedded in $X_0$.

An operator $A_0:X_1\to X_0$ is said to belong to the class $\mathcal{M}\mathcal{R}_p(X_1,X_0)$ if for each $f\in L_p(\mathbb{R}_+;X_0)$ there exists a unique solution
$$u\in H_p^1(\mathbb{R}_+;X_0)\cap L_p(\mathbb{R}_+;X_1)$$
of the problem
$$\dot{u}+A_0 u=f,\ t>0,\quad u(0)=0.$$
In other words, $A_0$ has maximal regularity of type $L_p$.

\section{Local Well-Posedness}\label{LWPsec}

The aim of this section is to solve the quasilinear evolution equation
    \begin{equation}\label{LWP1}
    \dot{u}+A(u)u=F_1(u)+F_2(u),\ t>0,\quad u(0)=u_1,
    \end{equation}
under the assumption that there exist two Banach spaces $X_0,X_1$, with dense embedding $X_1\hookrightarrow X_0$ such that the nonlinear mappings $(A,F_1)$ satisfy
    \begin{equation}\label{LWP2}
    (A,F_1)\in C^{1-}(V_\mu;\calB(X_1,X_0)\times X_0),
    \end{equation}
where $V_\mu\subset (X_0,X_1)_{\mu-1/p,p}=X_{\gamma,\mu}$ is open and nonempty for some $\mu\in (1/p,1]$. Furthermore, let $\beta\in (\mu-1/p,1)$ and $F_2:V_\mu\cap X_\beta\to X_0$ where $X_\beta=(X_0,X_1)_{\beta,p}$.

{\em Suppose that there exist numbers $m\in\mathbb{N}$, $\rho_j\ge0$ and $\beta_j\in [\mu-1/p,\beta]$ with
\begin{equation}\label{LWP:AssF_20}
\frac{\rho_j(\beta-\mu+1/p)+\beta_j-\mu+1/p}{1-\mu+1/p}<1,
\end{equation}
for all $j\in\{1,\ldots,m\}$, such that for each $u_*\in V_\mu$ and $R>0$ with $\bar{B}_R^{X_{\gamma,\mu}}(u_*)\subset V_\mu$ there exists $C_R>0$ such that the estimate
\begin{equation}\label{LWP:AssF_2}
|F_2(u_1)-F_2(u_2)|_{X_0}\le C_R\sum_{j=1}^m(1+|u_1|_{X_{\beta}}^{\rho_j}+|u_2|_{X_{\beta}}^{\rho_j})
|u_1-u_2|_{X_{\beta_j}},
\end{equation}
holds for all $u_1,u_2\in \bar{B}_R^{X_{\gamma,\mu}}(u_*)\cap X_\beta$, where $X_{\beta_j}=(X_0,X_1)_{\beta_j,p}$.}

The main result of this section reads as follows.
\begin{thm}\label{LWPthm}
Let $p\in (1,\infty)$ and suppose that $(A,F_1)$ and $F_2$ satisfy \eqref{LWP2}-\eqref{LWP:AssF_2} for some $\mu\in (1/p,1]$. Let $u_0\in V_\mu$ and assume that $A(u_0)\in\calM\calR_{p}(X_1,X_0)$. Then there exist $T=T(u_0)>0$ and $\ep=\ep(u_0)>0$, such that $\bar{B}_{\ep}^{X_{\gamma,\mu}}(u_0)\subset V_\mu$ and such that problem \eqref{LWP1} has a unique solution
    $$u(\cdot,u_1)\in H_{p,\mu}^1(0,T;X_0)\cap L_{p,\mu}(0,T;X_1)\cap C([0,T];V_\mu),$$
on $[0,T]$, for any initial value $u_1\in \bar{B}_\ep^{X_{\gamma,\mu}}(u_0)$. Furthermore there exists a constant $c=c(u_0)>0$ such that for all $u_1,u_2\in \bar{B}_\ep^{X_{\gamma,\mu}}(u_0)$ the estimate
    $$||u(\cdot,u_1)-u(\cdot,u_2)||_{\E_{1,\mu}(0,T)}\le c|u_1-u_2|_{X_{\gamma,\mu}}$$
is valid.
\end{thm}
\bpr
Since $u_0\in V_\mu$ and by \eqref{LWP2}, there exists $\ep_0>0$ and a constant $L>0$ such that $\bar{B}_{\ep_0}^{X_{\gamma,\mu}}(u_0)\subset V_\mu$ and
    \begin{equation}\label{LWP3}
    |A(w_1)v-A(w_2)v|_{X_0}\le L|w_1-w_2|_{X_{\gamma,\mu}}|v|_{X_1},
    \end{equation}
as well as
    \begin{equation}\label{LWP4}
    |F_1(w_1)-F_1(w_2)|_{X_0}\le L|w_1-w_2|_{X_{\gamma,\mu}},
    \end{equation}
hold for all $w_1,w_2\in \bar{B}_{\ep_0}^{X_{\gamma,\mu}}(u_0)$, $v\in X_1$. Moreover, from \eqref{LWP:AssF_2} we obtain the estimate
\begin{equation}\label{LWP:F_21}
|F_2(w_1)-F_2(w_2)|_{X_0}\le C_{\varepsilon_0}\sum_{j=1}^m(1+|w_1|_{X_{\beta}}^{\rho_j}+|w_2|_{X_{\beta}}^{\rho_j})
|w_1-w_2|_{X_{\beta_j}}
\end{equation}
for all $w_1,w_2\in \bar{B}_{\ep_0}^{X_{\gamma,\mu}}(u_0)\cap X_\beta$, where $\beta$ and $(\rho_j,\beta_j)$ satisfy \eqref{LWP:AssF_20}.

By \cite[Theorem 3.2]{PrSi04} we may introduce a reference function $u_0^*\in \E_{1,\mu}(0,T)$ as the solution of the linear problem
    $$\dot{w}+A(u_0)w=0,\quad w(0)=u_0.$$
Define a ball in $\E_{1,\mu}(0,T)$ by
    $$\B_{r,T,u_1}:=\{v\in \E_{1,\mu}(0,T):v|_{t=0}=u_1\ \mbox{and}\ ||v-u_0^*||_{\E_{1,\mu}(0,T)}\le r\},\quad 0<r\le1.$$
Let $u_1\in \bar{B}_{\ep}^{X_{\gamma,\mu}}(u_0)$ with $\ep\in (0,\ep_0]$. We will show that for all $v\in\B_{r,T,u_1}$ it holds that $v(t)\in \bar{B}_{\ep_0}^{X_{\gamma,\mu}}(u_0)$ for all $t\in [0,T]$, provided that $r,T,\ep>0$ are sufficiently small. For this purpose we define $u_1^*\in\E_{1,\mu}(0,T)$ as the unique solution of
    $$\dot{w}+A(u_0)w=0,\quad w(0)=u_1.$$
Given $v\in \B_{r,T,u_1}$ we estimate as follows.
    \begin{align}\label{LWP4a}
      ||v-u_0||_{\infty,X_{\gamma,\mu}}\le||v-u_1^*||_{\infty,X_{\gamma,\mu}}+||u_1^*-u_0^*||_{\infty,X_{\gamma,\mu}}
        +||u_0^*-u_0||_{\infty,X_{\gamma,\mu}}.
    \end{align}
Since $u_0$ is fixed, there exists $T_0=T_0(u_0)>0$ such that $\sup_{t\in[0,T_0]}|u_0^*(t)-u_0|_{X_{\gamma,\mu}}\le\ep_0/3$. Observe that $(v-u_1^*)|_{t=0}=0$, hence
    \begin{equation}\label{eq:tracethm}
    ||v-u_1^*||_{\infty,X_{\gamma,\mu}}\le C_1||v-u_1^*||_{\E_{1,\mu}(0,T)}
    \end{equation}
and the constant $C_1>0$ does not depend on $T$, which can be seen as follows. Thanks to \cite[Lemma 2.5]{MeySchn12} there exists a linear and bounded extension operator $\mathcal{E}_T^0$ from $_0\mathbb{E}_{1,\mu}(0,T)$ to $_0\mathbb{E}_{1,\mu}(\mathbb{R}_+)$, whose norm does not depend on $T$. This fact together with \cite[Proposition 3.1]{PrSi04} yields
$$\|u\|_{BUC(0,T;X_{\gamma,\mu})}\le \|\mathcal{E}_T^0 u\|_{BUC(\mathbb{R}_+;X_{\gamma,\mu})}\le C\|\mathcal{E}_T^0 u\|_{_0\mathbb{E}_{1,\mu}(\mathbb{R}_+)}\le C\|u\|_{_0\mathbb{E}_{1,\mu}(0,T)},$$
for each $u\in\,_0\mathbb{E}_{1,\mu}(0,T)$.

Making use of \eqref{eq:tracethm} we obtain
    \begin{align*}
    ||v-u_1^*||_{\infty,X_{\gamma,\mu}}&\le C_1||v-u_1^*||_{\E_{1,\mu}(0,T)}\le C_1(||v-u_0^*||_{\E_{1,\mu}(0,T)}+||u_0^*-u_1^*||_{\E_{1,\mu}(0,T)})\\
        &\le C_1(r+||u_0^*-u_1^*||_{\E_{1,\mu}(0,T)}),
    \end{align*}
and \eqref{LWP4a} yields the estimate
    \begin{multline*}
    ||v-u_0||_{\infty,X_{\gamma,\mu}}\\
        \le C_1(r+||u_0^*-u_1^*||_{\E_{1,\mu}(0,T)})+||u_0^*-u_1^*||_{\infty,X_{\gamma,\mu}}+||u_0^*-u_0||_{\infty,X_{\gamma,\mu}}.
    \end{multline*}
Since, by assumption, the semigroup $e^{-A(u_0)t}$ is exponentially stable it follows that
    \begin{equation}\label{LWP4b}
    ||u_0^*-u_1^*||_{\infty,X_{\gamma,\mu}}+C_1||u_0^*-u_1^*||_{\E_{1,\mu}(0,T)}\le C_\gamma|u_0-u_1|_{X_{\gamma,\mu}},
    \end{equation}
with a constant $C_\gamma>0$ which does not depend on $T$. Choosing $\ep\le\ep_0/(3C_\gamma)$ and $r\le\ep_0/(3C_1)$, we finally obtain
    \begin{equation}\label{LWP4c}
    ||v-u_0||_{\infty,X_{\gamma,\mu}}\le C_1r+C_\gamma\ep+||u_0^*-u_0||_{\infty,X_{\gamma,\mu}}\le\ep_0.
    \end{equation}
Throughout the remainder of this proof we will assume that $u_1\in B_{\ep}^{X_{\gamma,\mu}}(u_0)$, $\ep\le\ep_0/(3C_\gamma)$, $T\in[0,T_0]$ and $r\le\ep_0/(3C_1)$. Observe that in particular we have $v(t)\in \bar{B}_{\ep_0}^{X_{\gamma,\mu}}(u_0)\cap X_1$ for a.e.\ $t\in [0,T]$, hence also $v(t)\in \bar{B}_{\ep_0}^{X_{\gamma,\mu}}(u_0)\cap X_\beta$ for a.e.\ $t\in [0,T]$.

Under these assumptions, we may define a mapping $\calT_{u_1}:\B_{r,T,u_1}\to \E_{1,\mu}(0,T)$ by means of $\calT_{u_1} v=u$, where $u$ is the unique solution of the linear problem
    \begin{equation*}
    \dot{u}+A(u_0)u=F_1(v)+F_2(v)+(A(u_0)-A(v))v,\ t>0,\quad u(0)=u_1.
    \end{equation*}
In order to apply the contraction mapping principle, we have to show $\calT_{u_1}\B_{r,T,u_1}\subset\B_{r,T,u_1}$ and that $\calT_{u_1}$ defines a strict contraction on $\B_{r,T,u_1}$, i.e. there exists $\kappa\in(0,1)$ such that
    $$||\calT_{u_1} v-\calT_{u_1}\bar{v}||_{\E_{1,\mu}(0,T)}\le \kappa||v-\bar{v}||_{\E_{1,\mu}(0,T)},$$
is valid for all $v,\bar{v}\in\B_{r,T,u_1}$. We will first take care of the self-mapping property. Note that for $v\in\B_{r,T,u_1}$ we have
    $$\calT_{u_1} v-u_0^*=\calT_{u_1}v-u_1^*+u_1^*-u_0^*.$$
Since $(\calT_{u_1}v-u_1^*)|_{t=0}=0$, maximal $L_{p,\mu}$-regularity of $A(u_0)$ yields
    $$
    ||\calT_{u_1}v-u_1^*||_{\E_{1,\mu}(0,T)}
        \le C_0||F_1(v)+F_2(v)+(A(u_0)-A(v))v)||_{\E_{0,\mu}(0,T)},
    $$
and $C_0>0$ does not depend on $T$. It has been shown in \cite{KPW10} that
    $$
    ||(A(u_0)-A(v))v||_{\E_{0,\mu}(0,T)}\le L(C_1r+C_\gamma\ep+||u_0^*-u_0||_{\infty,X_{\gamma,\mu}})(r+||u_0^*||_{\E_{1,\mu}(0,T)})
    $$
and
    $$
    ||F_1(v)||_{\E_{0,\mu}(0,T)}\le
    \sigma(T)\left[L(C_1r+C_\gamma\ep+||u_0^*-u_0||_{\infty,X_{\gamma,\mu}})+|F_1(u_0)|_{X_0}\right],
    $$
with $\sigma(T):=\frac{1}{(1+(1-\mu)p)^{1/p}}T^{1/p+1-\mu}$. Since
    $$||u_0^*-u_0||_{\infty,X_{\gamma,\mu}},||u_0^*||_{\E_{1,\mu}(0,T)}\to 0\ \mbox{as}\ T\to 0_+,$$
this yields
    $$||\calT_{u_1}v-u_0^*||_{\E_{1,\mu}(0,T)}\le ||u_1^*-u_0^*||_{\E_{1,\mu}(0,T)}+\|F_2(v)\|_{\mathbb{E}_{0,\mu}(0,T)}+r/3,$$
provided $r>0,T>0,\ep>0$ are chosen properly. By \eqref{LWP4b} we obtain, in addition,
    $$||u_1^*-u_0^*||_{\E_{1,\mu}(0,T)}\le (C_\gamma/C_1)|u_1-u_0|_{X_{\ga,\mu}}\le r/3,$$
with a probably smaller $\ep>0$. Thus, it remains to estimate $F_2(v)$ in $\mathbb{E}_{0,\mu}(0,T)$. First we use the estimate
$$\|F_2(v)\|_{\mathbb{E}_{0,\mu}(0,T)}\le \|F_2(v)-F_2(u_0^*)\|_{\mathbb{E}_{0,\mu}(0,T)}+\|F_2(u_0^*)\|_{\mathbb{E}_{0,\mu}(0,T)}.$$
Then, by \eqref{LWP:F_21}, we obtain
\begin{multline}\label{LWP:EstF_2}
\|F_2(v)-F_2(u_0^*)\|_{\mathbb{E}_{0,\mu}(0,T)}\\
\le C_{\varepsilon_0}\sum_{j=1}^m\left(\int_0^T
(1+|v(t)|_{X_{\beta}}^{\rho_j}+|u_0^*(t)|_{X_{\beta}}^{\rho_j})^p
|v(t)-u_0^*(t)|_{X_{\beta_j}}^p
t^{(1-\mu)p}dt\right)^{1/p}
\end{multline}
To shorten the notation, let $w_1\in\{v,u_0^*\}$ and $w_2=v-u_0^*$. Observe that for all $w\in X_1$ we have the interpolation inequalities
\begin{equation}\label{eq:intineq}
|w|_{X_{\beta}}\le c|w|_{X_{\gamma,\mu}}^{1-\alpha}|w|_{X_1}^{\alpha}\quad\text{and}\quad
|w|_{X_{\beta_j}}\le c|w|_{X_{\gamma,\mu}}^{1-\alpha_j}|w|_{X_1}^{\alpha_j},
\end{equation}
where $\alpha(1-\mu+1/p)=(\beta-\mu+1/p)$ and $\alpha_j(1-\mu+1/p)=(\beta_j-\mu+1/p)$. It follows from \eqref{LWP:AssF_20} that $\rho_j\alpha+\alpha_j<1$ for all $j$. Making use of \eqref{eq:intineq} we obtain in a first step
\begin{multline*}
\left(\int_0^T |w_1(t)|_{X_{\beta}}^{\rho_jp}|w_2(t)|_{X_{\beta_j}}^{p}t^{(1-\mu)p}dt\right)^{1/p}\\
\le c\|w_1\|_{\infty,X_{\gamma,\mu}}^{\rho_j(1-\alpha)}
\|w_2\|_{\infty,X_{\gamma,\mu}}^{1-\alpha_j}\left(\int_0^T |w_1(t)|_{X_1}^{\rho_j\alpha p}|w_2(t)|_{X_1}^{\alpha_jp}t^{(1-\mu)p}dt\right)^{1/p}.
\end{multline*}
Appling the identity
$$t^{(1-\mu)p}=t^{(1-\mu)\rho_j\alpha p}t^{(1-\mu)\alpha_j p} t^{(1-\mu)(1-\rho_j\alpha-\alpha_j)p}$$
and H\"{o}lder's inequality twice yields
$$
\left(\int_0^T |w_1(t)|_{X_1}^{\rho_j\alpha p}|w_2(t)|_{X_1}^{\alpha_jp}t^{(1-\mu)p}dt\right)^{1/p}
\le T^{\kappa_{1,j}}\|w_1\|_{\mathbb{E}_{1,\mu}(0,T)}^{\rho_j\alpha}
\|w_2\|_{\mathbb{E}_{1,\mu}(0,T)}^{\alpha_j},
$$
with $\kappa_{1,j}:=(1-\mu+1/p)(1-\rho_j\alpha-\alpha_j)$. In a similar way we obtain
$$\left(\int_0^T|w_2(t)|_{\beta_j}^p t^{(1-\mu)p} dt\right)^{1/p}\le
cT^{\kappa_{2,j}}\|w_2\|_{\infty,X_{\gamma,\mu}}^{1-\alpha_j}
\|w_2\|_{\mathbb{E}_{1,\mu}(0,T)}^{\alpha_j},$$
with $\kappa_{2,j}:=(1-\mu+1/p)(1-\alpha_j)$.
Since $\|F_2(u_0^*)\|_{\mathbb{E}_{0,\mu}(0,T)}\to 0$ as $T\to 0_+$, $\|v-u_0^*\|_{\mathbb{E}_{1,\mu}(0,T)}\le r$,
$$\|v\|_{\infty,X_{\gamma,\mu}}\le \|v-u_0^*\|_{\infty,X_{\gamma,\mu}}+\|u_0^*\|_{\infty,X_{\gamma,\mu}}\le 4\varepsilon_0/3+\|u_0^*\|_{\infty,X_{\gamma,\mu}},$$
and
$$\|v\|_{\mathbb{E}_{1,\mu}(0,T)}\le \|v-u_0^*\|_{\mathbb{E}_{1,\mu}(0,T)}+\|u_0^*\|_{\mathbb{E}_{1,\mu}(0,T)}\le r+\|u_0^*\|_{\mathbb{E}_{1,\mu}(0,T)},$$
it follows that $\|F_2(v)\|_{\mathbb{E}_{0,\mu}(0,T)}\le r/3$ provided that $T>0$ is sufficiently small. This proves the self-mapping property of $\calT_{u_1}$.

Returning to the contraction property, let $u_1,u_2\in \bar{B}_{\ep}^{X_{\ga,\mu}}(u_0)$, $u_2^*(t):=e^{-A(u_0)t}u_2$ and $v_1\in \B_{r,T,u_1}$, $v_2\in \B_{r,T,u_2}$ be given. Then, by maximal $L_{p,\mu}$-regularity  of $A(u_0)$, we have
    \begin{align}\label{LWP5}
    \begin{split}
    ||\calT_{u_1}v_1-\calT_{u_2}v_2||_{\E_{1,\mu}(0,T)}&\le ||u_1^*-u_2^*||_{\E_{1,\mu}(0,T)}+C_0||F_1(v_1)-F_1(v_2)||_{\E_{0,\mu}(0,T)}\\
        &+C_0||(A(v_1)-A(u_0))(v_1-v_2)||_{\E_{0,\mu}(0,T)}\\
        &+C_0||(A(v_1)-A(v_2))v_2||_{\E_{0,\mu}(0,T)}\\
        &+C_0||F_2(v_1)-F_2(v_2)||_{\E_{0,\mu}(0,T)},
    \end{split}
    \end{align}
with a constant $C_0>0$ being independent of $T>0$. For the first term we can make use of \eqref{LWP4b} where $u_0$ and $u_0^*$ have to be replaced by $u_2$ and $u_2^*$, respectively. The second term can be treated as follows. By \eqref{LWP4}, we obtain
    $$||F_1(v_1)-F_1(v_2)||_{\E_{0,\mu}(0,T)}\le \sigma(T)L||v_1-v_2||_{\infty,X_{\ga,\mu}},$$
while \eqref{LWP4b} and the trace theorem \eqref{eq:tracethm} imply
    \begin{align}\label{LWP6}
    \begin{split}
    ||v_1-v_2||_{\infty,X_{\ga,\mu}}&\le ||v_1-v_2-(u_1^*-u_2^*)||_{\infty,X_{\ga,\mu}}+||u_1^*-u_2^*||_{\infty,X_{\ga,\mu}}\\
    &\le C_1||v_1-v_2-(u_1^*-u_2^*)||_{\E_{1,\mu}(0,T)}+C_\gamma|u_1-u_2|_{X_{\ga,\mu}}\\
    &\le C_1||v_1-v_2||_{\E_{1,\mu}(0,T)}+C_\gamma(1+C_1)|u_1-u_2|_{X_{\ga,\mu}}.
    \end{split}
    \end{align}
This yields
    $$||F_1(v_1)-F_1(v_2)||_{\E_{0,\mu}(0,T)}\le \sigma(T)L\left(C_1||v_1-v_2||_{\E_{1,\mu}(0,T)}+C_\gamma(1+C_1)|u_1-u_2|_{X_{\ga,\mu}}\right).$$
For the terms involving the operator $A$ in \eqref{LWP5} we make use of \eqref{LWP3} which results in
    \begin{multline*}
    ||(A(v_1)-A(u_0))(v_1-v_2)||_{\E_{0,\mu}(0,T)}+||(A(v_1)-A(v_2))v_2||_{\E_{0,\mu}(0,T)}\\
        \le L(||v_1-u_0||_{\infty,X_{\ga,\mu}}||v_1-v_2||_{\E_{1,\mu}(0,T)}+||v_1-v_2||_{\infty,X_{\ga,\mu}}||v_2||_{\E_{1,\mu}(0,T)}).
    \end{multline*}
Finally, making use of \eqref{LWP:AssF_2} and mimicking the estimates for $F_2$ from above, we obtain
\begin{align*}
||F_2(v_1)&-F_2(v_2)||_{\E_{0,\mu}(0,T)}\\
&\le C\sum_{j=1}^m\eta_j(T)\|v_1-v_2\|_{\infty,X_{\gamma,\mu}}^{1-\alpha_j}
\|v_1-v_2\|_{\mathbb{E}_{1,\mu}(0,T)}^{\alpha_j}\times\\
&\hspace{1cm}\times\Big(1+\|v_1\|_{\infty,X_{\gamma,\mu}}^{\rho_j(1-\alpha)}
\|v_1\|_{\mathbb{E}_{1,\mu}(0,T)}^{\rho_j\alpha}+
\|v_2\|_{\infty,X_{\gamma,\mu}}^{\rho_j(1-\alpha)}
\|v_2\|_{\mathbb{E}_{1,\mu}(0,T)}^{\rho_j\alpha}\Big)
\end{align*}
with $\eta_j(T)\to 0_+$ as $T\to 0_+$ for all $j\in\{1,\ldots,m\}$.

Observe that by Young's inequality and \eqref{LWP6} we have
\begin{align*}
\|v_1-v_2\|_{\infty,X_{\gamma,\mu}}^{1-\alpha_j}
\|v_1-v_2\|_{\mathbb{E}_{1,\mu}(0,T)}^{\alpha_j}&\le(1-\alpha_j)
\|v_1-v_2\|_{\infty,X_{\gamma,\mu}}+
\alpha_j\|v_1-v_2\|_{\mathbb{E}_{1,\mu}(0,T)}\\
&\le C\left(|u_1-u_2|_{X_{\gamma,\mu}}+\|v_1-v_2\|_{\mathbb{E}_{1,\mu}(0,T)}\right).
\end{align*}
By \eqref{LWP4c}, the term $||v_1-u_0||_{\infty,X_{\ga,\mu}}$ can be made as small as we wish by decreasing $r>0,T>0$ and $\ep>0$. Furthermore, for $j\in\{1,2\}$ we have
    $$||v_j||_{\infty,X_{\gamma,\mu}}\le ||v_j-u_0||_{\infty,X_{\gamma,\mu}}+|u_0|_{X_{\gamma,\mu}}\le \varepsilon_0+|u_0|_{X_{\gamma,\mu}}$$
and
    $$||v_j||_{\E_{1,\mu}(0,T)}\le ||v_j-u_0^*||_{\E_{1,\mu}(0,T)}+||u_0^*||_{\E_{1,\mu}(0,T)}\le r+||u_0^*||_{\E_{1,\mu}(0,T)},$$
hence $||v_j||_{\E_{1,\mu}(0,T)}$ is small, provided $r>0$ and $T>0$ are small enough. In summary, making use of \eqref{LWP6} and choosing $r>0, T>0$ and $\ep>0$ sufficiently small, we obtain a constant $c=c(u_0)>0$ such that the estimate
    \begin{equation}\label{LWP7}
    ||\calT_{u_1}v_1-\calT_{u_2}v_2||_{\E_{1,\mu}(0,T)}\le\frac{1}{2}||v_1-v_2||_{\E_{1,\mu}(0,T)}+c|u_1-u_2|_{X_{\ga,\mu}},
    \end{equation}
is valid for all $u_1,u_2\in \bar{B}_{\ep}^{X_{\ga,\mu}}(u_0)$ and $v_1\in \B_{r,T,u_1}$, $v_2\in \B_{r,T,u_2}$. In the very special case $u_1=u_2$, \eqref{LWP7} yields the contraction mapping property of $\calT_{u_1}$ on $\B_{r,T,u_1}$. Now we are in a position to apply Banach's fixed point theorem to obtain a unique fixed point $\tilde{u}\in\B_{r,T,u_1}$ of $\calT_{u_1}$, i.e. $\calT_{u_1}\tilde{u}=\tilde{u}$. Therefore $\tilde{u}\in\B_{r,T,u_1}$ is the unique local solution to \eqref{LWP1}. Furthermore, if $u(t,u_1)$ and $u(t,u_2)$ denote the solutions of \eqref{LWP1} with initial values $u_1,u_2\in \bar{B}_{\ep}^{X_{\ga,\mu}}(u_0)$, respectively, the last assertion of the theorem follows from \eqref{LWP7}.

\epr

\noindent
The next result provides information about the continuation of local solutions.
\begin{cor}\label{LWPcor1}
Let the assumptions of Theorem \ref{LWPthm} be satisfied and assume that $A(v)\in \calM\calR_{p}(X_1,X_0)$ for all $v\in V_\mu$. Then the solution $u(t)$ of \eqref{LWP1} with initial value $u_0\in V_\mu$ has a maximal interval of existence $J(u_0)=[0,t^+(u_0))$.

The mapping $[u_0\mapsto t^+(u_0)]:V_\mu\to (0,\infty)$ is lower-semicontinuous.
\end{cor}
\bpr
Given $u_0\in V_\mu\subset X_{\gamma,\mu}$, Theorem \ref{LWPthm} yields some $T_1>0$ and a unique solution $\bar{u}\in \E_{1,\mu}(0,T_1)\cap C([0,T_1];V_\mu)$ of \eqref{LWP1}. Next, we solve \eqref{LWP1} with initial value $\bar{u}(T_1)\in V_\mu$ to obtain a unique solution $\tilde{u}\in \E_{1,\mu}(0,T_2)\cap C([0,T_2];V_\mu)$ for some $T_2\in (0,\infty)$.
Let
    $$u(t):=\begin{cases}
            \bar{u}(t),&\quad t\in[0,T_1],\\
            \tilde{u}(t-T_1),&\quad t\in[T_1,T_1+T_2].
            \end{cases}
    $$
Then $u\in \E_{1,\mu}(0,T_1+T_2)\cap C([0,T_1+T_2];V_\mu)$, provided that
$$\int_{T_1}^{T_1+T_2}|\tilde{u}(t-T_1)|_1^p\  t^{(1-\mu)p}dt+\int_{T_1}^{T_1+T_2}|\dot{\tilde{u}}(t-T_1)|_0^p\ t^{(1-\mu)p}dt<\infty.$$
To show this property, it is sufficient to prove that $\tilde{u}\in \mathbb{E}_{1}(0,T_2)$. For that purpose we solve \eqref{LWP1} with initial value $\bar{u}(T_1-\delta)$ and some sufficiently small $\delta\in (0,T_2)$. This yields a unique solution $\hat{u}\in \E_{1,\mu}(0,T_2)\cap C([0,T_2];V_\mu)$ by Theorem \ref{LWPthm}, since the existence time $T_2>0$ is locally uniform.

Let us show that $\hat{u}(t)=\bar{u}(t+T_1-\delta)$ for $t\in [0,\delta]$. To this end we define $v(t):=\bar{u}(t+T_1-\delta)$ for $t\in [0,\delta]$. Then $v\in \E_{1,\mu}(0,\delta)$ and $v$ solves \eqref{LWP1} with $v(0)=\bar{u}(T_1-\delta)=\hat{u}(0)$. By uniqueness of the solution, it follows that $\hat{u}(t)=v(t)$ for each $t\in [0,\delta]$.

Now we show that $\tilde{u}(t)=\hat{u}(t+\delta)$ for $t\in [0,T_2-\delta]$. Let $w(t):=\hat{u}(t+\delta)$ for $t\in [0,T_2-\delta]$. Then $w\in \E_{1,\mu}(0,T_2-\delta)$ and $w$ solves \eqref{LWP1} with $w(0)=\hat{u}(\delta)=\bar{u}(T_1)=\tilde{u}(0)$. Therefore, by uniqueness, we have $w(t)=\tilde{u}(t)$ for each $t\in [0,T_2-\delta]$.

We claim that this already yields $\tilde{u}\in \mathbb{E}_{1}(0,T_2)$. Indeed, it holds that
\begin{align*}
\int_0^{T_2-\delta}|\tilde{u}(t)|_1^p dt&+\int_0^{T_2-\delta}|\dot{\tilde{u}}(t)|_0^pdt=
\int_0^{T_2-\delta}|\hat{u}(t+\delta)|_1^p dt+\int_0^{T_2-\delta}|\dot{\hat{u}}(t+\delta)|_0^pdt\\
&=\int_{\delta}^{T_2}|\hat{u}(\tau)|_1^p d\tau+\int_{\delta}^{T_2}|\dot{\hat{u}}(\tau)|_0^pd\tau\\
&\le \frac{1}{\delta^{(1-\mu)p}}\left(\int_{\delta}^{T_2}|\hat{u}(\tau)|_1^p \tau^{(1-\mu)p}d\tau+\int_{\delta}^{T_2}|\dot{\hat{u}}(\tau)|_0^p\tau^{(1-\mu)p}d\tau\right)\\
&\le \frac{1}{\delta^{(1-\mu)p}}\|\hat{u}\|_{\mathbb{E}_{1,\mu}(0,T_2)}^p.
\end{align*}
Additionally, we know $\mathbb{E}_{1,\mu}(T_2-\delta,T_2)\hookrightarrow \mathbb{E}_{1}(T_2-\delta,T_2)$, and hence it follows that $\tilde{u}\in \mathbb{E}_{1}(0,T_2)$. This in turn yields that $u\in \E_{1,\mu}(0,T_1+T_2)\cap C([0,T_1+T_2];V_\mu)$ is the unique solution of \eqref{LWP1} on the interval $[0,T_1+T_2]$. Inductively this yields a maximal interval of existence $J(u_0):=[0,t^+(u_0))$, which is of course half open, since otherwise we could continue the solution beyond $t^+(u_0)$ with initial value $u(t^+(u_0))$.

For the proof of the last assertion of the corollary one may follow the proof of \cite[Theorem 5.1]{CleSi01} (see also \cite[Theorem 7.2]{Ama88}). We refrain from repeating the arguments.

\epr

\begin{rem}
Let $J=[0,T_0]$ a compact interval and denote by $\calM\calR_p(J;X_1,X_0)$ the class of all linear operators $A_0:X_1\to X_0$ such that for all $f\in L_p(J;X_0)$ there exists a unique solution $u\in H_p^1(J;X_0)\cap L_p(J;X_1)$ of
    $$\dot{u}+A_0u=f,\ t\in (0,T_0],\quad u(0)=0.$$
It is well-known that this property does not depend on $T_0\in (0,\infty)$, and that there exists a number $\kappa>0$ such that the implication
    $$A_0\in \calM\calR_p(J;X_1,X_0)\Rightarrow A_0+\kappa I\in \calM\calR_p(X_1,X_0)$$
holds, see e.g.\ Pr\"{u}ss \cite{JanBari}. In this sense the assumption $A(u_0)\in\calM\calR_p(X_1,X_0)$ in Theorem \ref{LWPthm} can be replaced by the somewhat weaker condition $A(u_0)\in\calM\calR_p(J;X_1,X_0)$, we simply have to add $\kappa u$ to both sides of \eqref{LWP1} for some $\kappa>0$.
\end{rem}

\section{Relative compactness of orbits}\label{sec:relcomp}

Let $u_0\in V_\mu$ for some $\mu\in (1/p,1)$ be given. Suppose that $(A,F_1)$ and $F_2$ satisfy \eqref{LWP2}-\eqref{LWP:AssF_2}, respectively, and $A(v)\in \calM\calR_p(J;X_1,X_0)$ for all $v\in V_\mu$, where $J=[0,T]$ or $J=\R_+$. In the sequel we assume that the unique solution of \eqref{LWP1} satisfies $u\in BC([\tau,t^+(u_0));V_\mu\cap X_{\gamma,\bar{\mu}})$ for some $\tau\in (0,t^+(u_0))$, $\bar{\mu}\in (\mu,1]$ and
    \begin{equation}\label{GWPdist}
    \dist(u(t),\partial V_\mu)\ge\eta>0
    \end{equation}
for all $t\in J(u_0)=[0,t^+(u_0))$. Suppose furthermore that
    \begin{equation}\label{GWPemb}
    X_{\gamma,\bar{\mu}}\compemb X_{\gamma,\mu},\quad \bar{\mu}\in (\mu,1].
    \end{equation}
It follows from the boundedness of $u(t)$ in $X_{\gamma,\bar{\mu}}$ that the set $\{u(t)\}_{t\in J(u_0)}\subset V_\mu$ is relatively compact in $X_{\gamma,\mu}$, provided $\bar{\mu}\in (\mu,1]$. By \eqref{GWPdist} it holds that $\calV:=\overline{\{u(t)\}}_{t\in J(u_0)}$ is a proper subset of $V_\mu$. Applying Theorem \ref{LWPthm} we find for each $v\in \calV$ numbers $\ep(v)>0$ and $\de(v)>0$ such that $B_{\ep(v)}^{X_{\gamma,\mu}}(v)\subset V_\mu$ and all solutions of \eqref{LWP1} which start in $B_{\ep(v)}^{X_{\gamma,\mu}}(v)$ have the common interval of existence $[0,\de(v)]$. Therefore the set
    $$\bigcup_{v\in \calV}B_{\ep(v)}^{X_{\gamma,\mu}}(v)$$
is an open covering of $\calV$ and by compactness of $\calV$ there exist $N\in\N$ and $v_k\in\calV$, $k=1,\ldots,N$, such that
    $$\calU:=\bigcup_{k=1}^N B_{\ep_k}^{X_{\gamma,\mu}}(v_k)\supset\calV=\overline{\{u(t)\}}_{t\in J(u_0)}\supset\{u(t)\}_{t\in J(u_0)},$$
where $\ep_k:=\ep(v_k)$, $k=1,\ldots,N$. To each of these balls corresponds an interval of existence $[0,\de_k]$, $\de_k>0$, $k=1,\ldots,N$.
Consider the problem
    \begin{equation}\label{GWP4}
    \dot{v}+A(v)v=F_1(v)+F_2(v),\ s>0,\quad v(0)=u(t),
    \end{equation}
where $t\in J(u_0)$ is fixed and let $\de:=\min\{\de_k,\ k=1,\ldots,N\}$. Since $u(t)\subset \calU,\ t\in J(u_0)$, the solution of \eqref{GWP4} exists at least on the interval $[0,\de]$. Assume that $t^+(u_0)<\infty$ and set $t=t^+(u_0)-\delta/2$ in \eqref{GWP4}. Then we obtain a unique solution $v(s)$ of \eqref{GWP4}  for $s\in[0,\delta]$, hence
$$w(\tau):=\begin{cases}
u(\tau),\ \tau\in[0,t^+(u_0)-\delta/2],\\
v(\tau-t^+(u_0)+\delta/2),\ \tau\in[t^+(u_0)-\delta/2,t^+(u_0)+\delta/2],
\end{cases}$$
is the unique solution of \eqref{LWP1} on the interval $[0,t^+(u_0)+\delta/2]$. This contradicts the maximality of $t^+(u_0)$, hence the solution exists globally.

By continuous dependence on the initial data, the solution operator $G_1:\mathcal{U}\to \E_{1,\mu}(0,\de)$, which assigns to each initial value $u_1\in \mathcal{U}$ a unique solution $v(\cdot,u_1)\in \E_{1,\mu}(0,\de)$, is continuous. Furthermore
    $$(\de/2)^{1-\mu}||v||_{\E_{1}(\de/2,\de)}\le ||v||_{\E_{1,\mu}(\de/2,\de)}\le ||v||_{\E_{1,\mu}(0,\de)},\ \mu\in (1/p,1),$$
hence the mapping $G_2:\E_{1,\mu}(0,\de)\to \E_{1}(\de/2,\de)$ with $v\mapsto v$ is continuous. Finally
    $$|v(\delta)|_{X_{\gamma}}\le||v||_{BUC((\de/2,\de);X_{\gamma})}\le C(\delta)||v||_{\E_1(\de/2,\de)},$$
which implies that the mapping $G_3:\E_1(\de/2,\de)\to X_{\gamma}$ with $v\mapsto v(\de)$ is continuous. This yields the continuity of the composition $G=G_3\circ G_2\circ G_1:\mathcal{U}\to X_{\gamma}$, whence $G(\{u(t)\}_{t\ge 0})=\{u(t+\de)\}_{t\ge 0}$ is relatively compact in $X_{\gamma}$, since the continuous image of a relatively compact set is relatively compact. Since the solution has relatively compact range in $X_\gamma$, it is an easy consequence that the $\om$-limit set
    \begin{equation}\label{omegalimset}
    \om(u_0):=\left\{v\in V_\mu\cap X_\gamma:\ \exists\ t_n\nearrow\infty\ \mbox{s.t.}\ u(t_n;u_0)\to v\ \mbox{in}\ X_\gamma\right\}
    \end{equation}
is nonempty, compact and connected. We summarize the preceding considerations in the following
\begin{thm}\label{GWPthm2}
Let $p\in (1,\infty)$ and let $J=[0,T]$ or $J=\R_+$. Suppose that $A(v)\in\calM\calR_{p}(J;X_1,X_0)$ for all $v\in V_\mu$ and let \eqref{LWP2}-\eqref{LWP:AssF_2} hold for some $\mu\in (1/p,1)$. Assume furthermore that \eqref{GWPemb} holds for some $\bar{\mu}\in (\mu,1]$ and that the solution $u(t)$ of \eqref{LWP1} satisfies
    $$u\in BC([\tau,t^+(u_0));V_\mu\cap X_{\gamma,\bar{\mu}})$$
for some $\tau\in (0,t^+(u_0))$ and $\bar{\mu}\in (\mu,1]$ such that
    $$\dist(u(t),\pa V_\mu)\ge \eta>0$$
for all $t\in J(u_0)$. Then the solution exists globally and for each $\de>0$, the orbit $\{u(t)\}_{t\ge \de}$ is relatively compact in $X_{\gamma}$. If in addition $u_0\in V_\mu\cap X_{\gamma}$, then $\{u(t)\}_{t\ge 0}$ is relatively compact in $X_{\gamma}$.

In either case, the $\omega$-limit set $\omega(u_0)$ given by \eqref{omegalimset} is nonempty, compact and connected.
\end{thm}

\section{Applications}\label{sec:appl}

\subsection{Reaction-Diffusion Systems}

Let $\Omega\subset\mathbb{R}^n$, $n\in\mathbb{N}$, be a bounded domain with boundary $\partial\Omega\in C^2$ and $E$ a finite dimensional real Hilbert space. Furthermore, let $U\subset E$ be open and suppose that $f\in C^{1-}(U;E)$ and $a\in C^{2-}(U;\mathcal{B}(E))$ are given. The general form of systems we consider here is given by
$$\partial_t u-\diver (a(u)\nabla u)=f(u),$$
where $u:(0,\infty)\times\Omega\to U$ and
$$\diver(a(u)\nabla u):=\sum_{j=1}^n \partial_j(a(u)\partial_j u).$$
This includes reaction-diffusion systems of Maxwell-Stefan type, see \cite{Both11}.
Computing $\diver(a(u)\nabla u)$, we obtain
$$\partial_t u-a(u)\Delta u=f(u)+\sum_{j=1}^n a'(u)\partial_ju\partial_j u.$$
This motivates to consider the problem
    \begin{align}
    \begin{split}\label{appl1}
    \partial_t u-a(u)\Delta u&=f(u)+\sum_{j=1}^nb(u)\partial_j u\partial_j u,\ t>0,\ x\in\Omega,\\
    \partial_\nu u&=0,\ t>0,\ x\in\pa\Omega,\\
    u(0)&=u_0,\ x\in\Omega,
    \end{split}
    \end{align}
where
$$(a,f)\in C^{1-}(U;\mathcal{B}(E)\times E)\ \text{and}\ b\in C^{1-}(U;\mathcal{B}(E,\mathcal{B}(E))).$$
The boundary condition $\partial_\nu u$ is defined by
$$\partial_\nu u:=(\nu\cdot\nabla)u.$$
Let us first rewrite \eqref{appl1} in the form \eqref{LWP1}. To this end we set $X_0=L_q(\Om;E)$,
    $$X_1=\{u\in H_q^2(\Om;E):\partial_\nu u|_{\pa\Om}=0\}.$$
In this situation, we have for $\mu\in (1/p,1]$
    \begin{multline*}
    X_{\ga,\mu}=(X_0,X_1)_{\mu-1/p,p}\\
                                        =\begin{cases}
                                        \{u\in B_{qp}^{2\mu-2/p}(\Om;E):\partial_\nu u|_{\pa\Om}=0\},&\ \mbox{if}\ 2\mu>1+2/p+1/q ,\\
                                        B_{qp}^{2\mu-2/p}(\Om;E),&\ \mbox{if}\ 2\mu<1+2/p+1/q.
                                        \end{cases}
    \end{multline*}
Let us assume that $2/p+n/q<2$, wherefore the embedding $B_{qp}^{2-2/p}(\Om;E)\hookrightarrow C(\overline{\Om};E)$ is at our disposal. In this case there exists $\mu_0\in (1/p,1)$ such that
    $$B_{qp}^{2-2/p}(\Om;E)\compemb B_{qp}^{2\mu-2/p}(\Om;E)\hookrightarrow C(\overline{\Om};E),\quad\mbox{if}\ \mu\in (\mu_0,1).$$
Indeed, the number $\mu_0\in (1/p,1)$ is given by
    $$\mu_0=\frac{1}{p}+\frac{n}{2q}.$$
Define
$$V_{\mu}:=\{v\in X_{\gamma,\mu}:v(\overline{\Omega})\subset U\},$$
which is an open set in $X_{\gamma,\mu}$, since $X_{\gamma,\mu}\hookrightarrow C(\overline{\Omega};E)$.

In the context of Maxwell-Stefan diffusion one typically has $\sigma(a(v(x)))\subset (0,\infty)$ for all $v\in V_\mu$ and $x\in\overline{\Omega}$. This leads to the assumption $\sigma (a(u))\subset (0,\infty)$ for each $u\in U$.
For $\mu\in (\mu_0,1]$, we define
    $$A(v)u(x):=a(v(x))\Delta u(x),\ x\in\Om,\ v\in V_\mu,\ u\in X_1,$$
    $$F_1(v)(x):=f(v(x)),\ x\in\Om,\ v\in V_\mu,$$
and
    $$F_2(v)(x):=\sum_{j=1}^nb(v(x))\partial_jv(x)\partial_jv(x),\ x\in\Omega,\ v\in X_\beta\cap V_\mu,$$
    where $\beta\in (\mu-1/p,1)$ will be chosen in a suitable way such that $F_2$ satisfies the assumptions \eqref{LWP:AssF_20} and \eqref{LWP:AssF_2}.
From the conditions on $(a,f)$ we infer
    $$(A,F_1)\in C^{1-}(V_\mu;\calB(X_1,X_0)\times X_0),\quad \mu\in (\mu_0,1],$$
and, for each $v\in V_\mu$, $x_0\in\overline{\Omega}$, the symbol of $\mathcal{A}(x_0,D):=a(v(x_0))\Delta$ is normally elliptic. Moreover, the pair $(\mathcal{A}(x_0,D),\partial_\nu)$ satisfies the Lopatinskii-Shapiro condition (see e.g.\ \cite[Section 2]{DHP2}) for each $v\in V_\mu$ and $x_0\in\partial\Omega$, since the matrix $a(v(x_0))$ is invertible and $\sigma(a(v(x_0)))\subset (0,\infty)$.

Let us now verify the condition \eqref{LWP:AssF_2}. We know that $X_\beta=(X_0,X_1)_{\beta,p}\subset B_{qp}^{2\beta}(\Omega;E)$. In the sequel, let
$$\beta\in\left(1/2+n/4q,1\right)\cap (\mu-1/p,1).$$
Note that the first interval is not empty, since $2/p+n/q<2$.
For this choice of $\beta$ we have the embedding $B_{qp}^{2\beta}(\Omega;E)\hookrightarrow W^1_{2q}({\Omega};E)$ at our disposal, since $2\beta-n/q>1-n/2q$, and it holds that $X_\beta\hookrightarrow X_{\gamma,\mu}$. For a given $u_*\in V_\mu$, choose $R>0$ such that $\bar{B}_R^{X_{\gamma,\mu}}(u_*)\subset V_\mu$.
Since $X_{\gamma,\mu}\hookrightarrow C(\bar{\Omega};E)$ and $b$ is locally Lipschitz continuous there exist constants $C_R,L_R>0$ such that $|b(v)|_\infty\le C_R$ and
$$|b(v)-b(\bar{v})|_\infty\le L_R|v-\bar{v}|_\infty\le CL_R|v-\bar{v}|_{X_{\gamma,\mu}},$$
for all $v,\bar{v}\in \bar{B}_R^{X_{\gamma,\mu}}(u_*)\cap X_\beta$.
Next, we write
\begin{multline*}
b(v)\partial_jv\partial_jv-b(\bar{v})\partial_j\bar{v}\partial_j\bar{v}=
b(v)(\partial_jv\partial_jv-\partial_j\bar{v}\partial_j\bar{v})+
(b(v)-b(\bar{v}))\partial_j\bar{v}\partial_j\bar{v}\\
=b(v)[(\partial_jv-\partial_j\bar{v})\partial_j v+(\partial_jv-\partial_j\bar{v})\partial_j \bar{v}]+
(b(v)-b(\bar{v}))\partial_j\bar{v}\partial_j\bar{v}.
\end{multline*}
Making use of H\"{o}lder's inequality and the fact that $X_\beta\hookrightarrow W_{2q}^1(\Omega;E)$ we obtain the estimates
$$
\left(\int_\Omega|b(v(x))(\partial_j v(x)-\partial_j \bar{v}(x))\partial_j v(x)|_E^qdx\right)^{1/q}
\le C_R|v|_{X_\beta}|v-\bar{v}|_{X_\beta},
$$
$$
\left(\int_\Omega|b(v(x))(\partial_j v(x)-\partial_j \bar{v}(x))\partial_j \bar{v}(x)|_E^qdx\right)^{1/q}
\le C_R|\bar{v}|_{X_\beta}|v-\bar{v}|_{X_\beta}.
$$
and
$$\left(\int_\Omega|(b(v(x))-b(\bar{v}(x)))\partial_j \bar{v}(x)\partial_j \bar{v}(x)|_E^{q}dx\right)^{1/q}
\le CL_R|\bar{v}|_{X_\beta}^2|v-\bar{v}|_{X_{\gamma,\mu}},$$
for each $j\in\{1,\ldots,n\}$. Therefore, $F_2$ satisfies the estimate \eqref{LWP:AssF_2} with $m=2$, $(\rho_1,\beta_1)=(1,\beta)$ and $(\rho_2,\beta_2)=(2,\mu-1/p)$. It remains to verify \eqref{LWP:AssF_20}, i.e.\ we have to check that
$$\beta<\frac{1}{2}\left(1+\mu-\frac{1}{p}\right).$$
Since
$$\max\left\{\frac{1}{2}+\frac{n}{4q},\mu-\frac{1}{p}\right\}<
\frac{1}{2}\left(1+\mu-\frac{1}{p}\right),$$
for $\mu\in(\mu_0,1)$, it is always possible to find such a number $\beta$.

We are now in a position to apply Theorem \ref{LWPthm} which yields the existence of a unique solution $u$ of \eqref{appl1} with a maximal interval of existence $J(u_0)$, provided $u_0\in X_{\gamma,\mu}$. Then we may apply Theorem \ref{GWPthm2} to the result
\begin{thm}\label{thm:appl1}
Let $n\in\N$, $p,q\in (1,\infty)$ with $2/p+n/q<2$, $\mu_0:=1/p+n/2q$, $\mu\in (\mu_0,1)$, $\Omega\subset\R^n$ a bounded domain with boundary $\pa\Om\in C^2$. Suppose that $E$ is a finite dimensional real Hilbert space and let $U\subset E$ be open. Assume that $(a,f)\in C^{1-}(U;\mathcal{B}(E,E)\times E)$ with $\sigma(a(u))\subset (0,\infty)$ for each $u\in U$ and $b\in C^{1-}(U;\mathcal{B}(E;\mathcal{B}(E)))$. Let $u_0\in B_{qp}^{2\mu-2/p}(\Om;E)$ with $u_0(\overline{\Omega})\subset U$ and $\partial_\nu u_0|_{\pa\Om}=0$ if $2\mu>1+2/p+1/q$. Then there exists a unique solution $u(t)$ of \eqref{appl1} with
$$u\in H_{p,\mu}^1(0,T;L_q(\Omega;E))\cap L_{p,\mu}(0,T;H_q^2(\Omega;E))\cap C([0,T];B_{qp}^{2\mu-2/p}(\Omega;E))$$
and $u([0,T]\times\overline{\Omega})\subset U$ for each $T\in (0,t^+(u_0))$, where $t^+(u_0)>0$ is the maximal time of existence.

If in addition the solution satisfies
    $$u\in BC\left([\tau,t^+(u_0));B_{qp}^{2\bar{\mu}-2/p}(\Om;E)\right),$$
for some $\tau\in (0,t^+(u_0))$, $\bar{\mu}\in (\mu,1]$ such that
    $$\dist(u(t,x),\partial U)\ge\eta>0$$
for all $t \in [0,t^+(u_0))$, $x\in\overline{\Omega}$ with some $\eta>0$, then $u(t)$ exists globally and the set $\{u(t)\}_{t\ge \tau}$ is relatively compact in $B_{qp}^{2-2/p}(\Om;E)$. Moreover, the $\om$-limit set
    $$\om(u_0):=\left\{v\in B_{qp}^{2-2/p}(\Om;E):\ \exists\ t_n\nearrow\infty\ \mbox{s.t.}\ u(t_n,u_0)\to v\ \mbox{in}\ B_{qp}^{2-2/p}(\Om;E)\right\}$$
is nonempty, connected and compact.
\end{thm}

%\begin{rem}
%For simplicity we supplied $\eqref{appl1}_1$ with a Dirichlet boundary condition. However, this boundary condition may be replaced by any other one (up to differential order one), as long as the Lopatinskii-Shapiro condition holds, which leads to maximal $L_p$-regularity (see e.g.\ \cite{DHP1,DHP2}).
%\end{rem}
\begin{rem}\
As a consequence, to obtain global existence for system \eqref{appl1} we need H\"{o}lder-bounds for $u$, but not for $\nabla u$. This coincides with the results in \cite{Ama89} for the global existence of solutions to quasilinear parabolic equations.
\end{rem}

\subsection{Surface diffusion flow - The graph case}\label{subsec:surfdiff}

Consider a family of hypersurfaces $\{\Gamma(t)\}_{t\ge 0}\subset\mathbb{R}^{n+1}$ given as a graph of a height function $h$ over a bounded domain $\Omega\subset\mathbb{R}^{n}$ with boundary $\partial\Omega\in C^4$. To be precise we have
$$\Gamma(t)=\{(x,x_{n+1})\in\Omega\times\mathbb{R}:x_{n+1}=h(t,x)\}.$$
Assume furthermore that the evolution of $\Gamma(t)$ is governed by the \emph{surface diffusion law}
\begin{equation}\label{eq:appl2}
V_{\Gamma(t)}=-\Delta_{\Gamma(t)}H_{\Gamma(t)},\ t\ge 0,
\end{equation}
where $V_{\Gamma(t)}$ is the normal velocity of $\Gamma(t)$, $H_{\Gamma(t)}=-{\diver}_{\Gamma(t)}\nu_{\Gamma(t)}$ is the mean curvature of $\Gamma(t)$ and the operators ${\diver}_{\Gamma(t)}$ and $\Delta_{\Gamma(t)}$ denote the surface divergence and the Laplace-Beltrami operator, respectively, acting on $\Gamma(t)$. Let us make the convention that the unit normal field $\nu_{\Gamma(t)}$ on $\Gamma(t)$ points from
$$\{(x,x_{n+1})\in\Omega\times\mathbb{R}:x_{n+1}<h(t,x)\}$$
to
$$\{(x,x_{n+1})\in\Omega\times\mathbb{R}:x_{n+1}>h(t,x)\}.$$
It is convenient to rewrite \eqref{eq:appl2} in terms of the height function $h$. To this end let $\beta:=1/\sqrt{1+|\nabla h|^2}$ and denote by $\delta^{ij}$ the Kronecker delta. Then we obtain from \cite[Section 2]{PrSi13} that
$$\Delta_{\Gamma}\varphi=(\delta^{kl}-\beta^2\partial_k h\partial_l h)(\partial_k\partial_l\varphi-
\beta^2\partial_k\partial_l h\,\partial_m h\,\partial_m\varphi),$$
for functions $\varphi$ which are smooth enough and
$$H_{\Gamma}=(\delta^{ij}-\beta^2\partial_i h\partial_j h)\beta\partial_i\partial_j h.$$
Note that for the sake of readability we make use of sum convention. Furthermore we have
$$\nu_\Gamma=\beta(-\nabla h,1)^{\sf T}\quad \text{and}\quad V_\Gamma=\partial_t h(e_{n+1}|\nu_\Gamma)=\beta\partial_t h.$$
In terms of boundary conditions for the height function we choose $h|_{\partial\Omega}=\partial_\nu h=0$ at $\partial\Omega$, where $\nu$ denotes the outer unit normal field on $\partial\Omega$.

Inserting the above expressions into \eqref{eq:appl2} yields the equation
\begin{multline}\label{eq:appl2aa}
\partial_t h+\sum_{|\sigma|=4}a_\sigma(\nabla h)D^\sigma h
=\sum_{|\sigma|=3,|\tau|=2}b_{\sigma\tau}(\nabla h)D^\sigma h D^\tau h\\
+\sum_{|\sigma|=|\tau|=|\chi|=2}c_{\sigma\tau\chi}(\nabla h)D^\sigma h D^\tau h D^\chi h,\quad t>0,\ x\in\Omega,
\end{multline}
supplemented with the boundary and initial conditions
\begin{alignat}{2}\label{eq:appl2aaa}
h|_{\partial\Omega}&=0,&&\quad t>0,\ x\in\partial\Omega\nn,\\
\partial_\nu h&=0,&&\quad t>0,\ x\in\partial\Omega,\\
h(0)&=h_0,&&\quad x\in\Omega\nn,
\end{alignat}
where
$D^\alpha:=(-i)^{|\alpha|}\partial_1^{\alpha_1}\cdots\partial_{n}^{\alpha_n}$ and $\alpha\in\mathbb{N}_0^n$ is a multiindex. Note that the coefficients $a_\sigma,b_{\sigma\tau},c_{\sigma\tau\chi}$ are smooth and that the leading coefficient $a_\sigma(\nabla h)$ is given by
\begin{equation}\label{eq:appl2a}
a_\sigma(\nabla h)=(\delta^{kl}-\beta^2\partial_k h\partial_l h)(\delta^{ij}-\beta^2\partial_ih\partial_j h),\ i,j,k,l\in\{1,\ldots,n\}.
\end{equation}
Let us reformulate problem \eqref{eq:appl2aa} in the form \eqref{LWP1}. To this end let $1<q<\infty$ and define $X_0:=L_q(\Omega)$,
$$X_1:=\{h\in H_q^4(\Omega):h=\partial_\nu h=0\ \text{at}\ \partial\Omega\},$$
as well as $X_{\gamma,\mu}:=(X_0,X_1)_{\mu-1/p,1}$ for $\mu\in (1/p,1]$, which yields $X_{\gamma,\mu}\subset B_{qp}^{4\mu-4/p}(\Omega)$. Let us assume that $4/p+n/q<3$, wherefore the embedding $B_{qp}^{4-4/p}(\Omega)\hookrightarrow C^1(\overline{\Omega})$ is at our disposal. This readily yields
$$B_{qp}^{4-4/p}(\Omega)\compemb B_{qp}^{4\mu-4/p}(\Omega)\hookrightarrow C^1(\overline{\Omega}),$$
provided that $\mu\in (\mu_0,1)$, where
$$\mu_0:=\frac{1}{p}+\frac{n}{4q}+\frac{1}{4}\in (1/p,1).$$
In the sequel we will always assume that $\mu\in(\mu_0,1)$. In this case one has the characterization
$$X_{\gamma,\mu}=\{h\in B_{qp}^{4\mu-4/p}(\Omega):h|_{\partial\Omega}=\partial_\nu h=0\ \text{at}\ \partial\Omega\}.$$
For $h\in X_{\gamma,\mu}$ and $u\in X_1$ we define $F_1(h)=0$,
$$A(h)u(x):=\sum_{|\sigma|=4}a_\sigma(\nabla h(x))D^\sigma u(x)$$
and for $h\in X_\beta=(X_0,X_1)_{\beta,p}$
\begin{multline*}
F_2(h)(x):=\sum_{|\sigma|=3,|\tau|=2}b_{\sigma\tau}(\nabla h(x))D^\sigma h(x) D^\tau h(x)\\
+\sum_{|\sigma|=|\tau|=|\chi|=2}c_{\sigma\tau\chi}(\nabla h(x))D^\sigma h(x) D^\tau h(x) D^\chi h(x),
\end{multline*}
where $\beta\in (\mu-1/p,1)$ will be chosen in a suitable way such that $F_2$ satisfies the assumptions \eqref{LWP:AssF_20} and \eqref{LWP:AssF_2}.
In the sequel, let
\begin{equation}\label{eq:appl2ab}
\beta\in (3/4+n/12q,1)\cap (\mu-1/p,1),
\end{equation}
which is possible, since $4/p+n/q<3$. This yields the validity of the embedding $X_\beta\hookrightarrow W_{3q/2}^3(\Omega)$, since $4\beta-n/q>2-2n/3q$. Furthermore we have $W_{3q/2}^3(\Omega)\hookrightarrow W_{3q}^2(\Omega)$, since $4/p+n/q<3$.

\noindent
Let $h,\bar{h}\in \bar{B}_R^{X_{\gamma,\mu}}(h_*)\cap X_\beta$ for some $h_*\in X_{\gamma,\mu}$ and observe that
\begin{multline*}
b_{\sigma\tau}(\nabla h)D^\sigma h D^\tau h-b_{\sigma\tau}(\nabla \bar{h})D^\sigma \bar{h} D^\tau \bar{h}\\
=b(\nabla h)(D^\sigma h D^\tau h-D^\sigma \bar{h} D^\tau \bar{h})+D^\sigma \bar{h} D^\tau \bar{h}(b(\nabla h)-b(\nabla \bar{h}))\\
=b(\nabla h)[D^\sigma h(D^\tau h-D^\tau\bar{h})+D^\tau\bar{h}(D^\sigma h-D^\sigma\bar{h})]+D^\sigma \bar{h} D^\tau \bar{h}(b(\nabla h)-b(\nabla \bar{h})),
\end{multline*}
hence H\"{o}lder's inequality yields the estimate
\begin{multline}\label{eq:appl3}
|b_{\sigma\tau}(\nabla h)D^\sigma h D^\tau h-b_{\sigma\tau}(\nabla \bar{h})D^\sigma \bar{h} D^\tau \bar{h}|_{L_q}\\
\le |b_{\sigma\tau}(\nabla h)|_\infty[|h|_{W_{3q/2}^3}|h-\bar{h}|_{W_{3q}^2}+|\bar{h}|_{W_{3q}^2}
|h-\bar{h}|_{W_{3q/2}^3}]\\
+|b(\nabla h)-b(\nabla\bar{h})|_\infty|\bar{h}|_{W_{3q/2}^3}|\bar{h}|_{W_{3q}^2}
\end{multline}
Next, let $\kappa:=1/2+n/6q$. Then, for $\varepsilon>0$ with $\kappa+\varepsilon<1$, we have
$$X_{\kappa+\varepsilon}=(X_0,X_1)_{\kappa+\varepsilon,p}\subset B_{qp}^{4(\kappa+\varepsilon)}(\Omega)\hookrightarrow
W_{3q}^2(\Omega),$$
since $4(\kappa+\varepsilon)-n/q>2-n/3q$.
For the first term on the right side of \eqref{eq:appl3} we choose $(\rho_1,\beta_1)=(1,\kappa+\varepsilon)$. Note that $\beta_1=\kappa+\varepsilon<\beta$ provided $\varepsilon>0$ is sufficiently small. Substituting $(\rho_1,\beta_1)$ into \eqref{LWP:AssF_20} yields the restriction
\begin{equation}\label{eq:appl4}
\beta<\mu-\frac{1}{p}+\frac{1}{2}-\frac{n}{6q}-\varepsilon.
\end{equation}
It is easy to check that one can always find a number $\beta$ satisfying \eqref{eq:appl2ab} and \eqref{eq:appl4}, whenever $\mu\in (\mu_0,1)$, $4/p+n/q<3$ and $\varepsilon>0$ is sufficiently small.

Concerning the second term on the right side of \eqref{eq:appl3} note that
\begin{equation}\label{eq:appl5}
|\bar{h}|_{W_{3q}^2(\Omega)}\le c|\bar{h}|_{X_{\kappa+\varepsilon}}\le c|\bar{h}|_{X_\beta}^\theta|\bar{h}|_{X_{\gamma,\mu}}^{1-\theta},
\end{equation}
by the reiteration theorem, where $\theta(\beta-\mu+1/p)=\kappa+\varepsilon-\mu+1/p$. We set $(\rho_2,\beta_2)=(\theta,\beta)$ for the estimate of the second term. Then \eqref{LWP:AssF_20} becomes \eqref{eq:appl4} again.

To estimate the third term, we observe that
$$|\bar{h}|_{W_{3q/2}^3(\Omega)}|\bar{h}|_{W_{3q}^2(\Omega)}\le
c|\bar{h}|_{X_{\gamma,\mu}}^{1-\theta}|\bar{h}|_{X_\beta}^{1+\theta}\le C_R|\bar{h}|_{X_\beta}^{1+\theta},$$
by \eqref{eq:appl5}. Therefore the choice $(\rho_3,\beta_3)=(1+\theta,\mu-1/p)$ leads again to \eqref{eq:appl4}.

We will now deal with the second sum in the definition of $F_2$. H\"{o}lder's inequality and \eqref{eq:appl5} yield the estimate
\begin{align*}
|c_{\sigma\tau\chi}(\nabla h)[D^\sigma h D^\tau h D^\chi h-D^\sigma \bar{h} D^\tau \bar{h} D^\chi \bar{h}]|_{L_q}&\le C_R(|h|_{W_{3q}^2}^2+|\bar{h}|_{W_{3q}^2}^2)|h-\bar{h}|_{W_{3q}^2}\\
&\le C_R(|h|_{X_\beta}^{2\theta}+|\bar{h}|_{X_\beta}^{2\theta})|h-\bar{h}|_{X_{\kappa+\varepsilon}}.
\end{align*}
Therefore we choose this time $(\rho_4,\beta_4)=(2\theta,\kappa+\varepsilon)$, which yields that \eqref{LWP:AssF_20} is equivalent to $\mu>\mu_0+3\varepsilon/2$ which can be achieved, provided that $\varepsilon>0$ is sufficiently small.

Finally, we have
\begin{align*}
|[c_{\sigma\tau\chi}(\nabla h)-c_{\sigma\tau\chi}(\nabla \bar{h})]D^\sigma \bar{h} D^\tau \bar{h} D^\chi \bar{h}|_{L_q}&\le |\bar{h}|_{W_{3q}^2}^3|c_{\sigma\tau\chi}(\nabla h)-c_{\sigma\tau\chi}(\nabla \bar{h})|_\infty\\
&\le C_R|\bar{h}|_{X_{\beta}}^{3\theta}|c_{\sigma\tau\chi}(\nabla h)-c_{\sigma\tau\chi}(\nabla \bar{h})|_\infty,
\end{align*}
where me made again use of H\"{o}lder's inequality and \eqref{eq:appl5}. This time we set $(\rho_5,\beta_5)=(3\theta,\mu-1/p)$ which implies again that \eqref{LWP:AssF_20} is equivalent to the condition $\mu>\mu_0+3\varepsilon/2$.

In summary it follows that $F_2$ satisfies \eqref{LWP:AssF_20} and \eqref{LWP:AssF_2} with $m=5$ and the above choices of $\beta$ and $(\rho_j,\beta_j)$, $j\in\{1,\ldots,5\}$, since the coefficients $b_{\sigma\tau}$ and $c_{\sigma\tau\chi}$ are locally Lipschitz continuous in $\nabla h$.

It remains to check the conditions on $A$. Clearly, the leading coefficient $a_\sigma$ is continuously differentiable with respect to $\nabla h$, wherefore $a_\sigma$ is locally Lipschitz continuous. This in turn implies that
$$A\in C^{1-}(X_{\gamma,\mu};\mathcal{B}(X_1,X_0)),$$
since $X_{\gamma,\mu}\hookrightarrow C^{1}(\overline{\Omega})$. Next we show that $A(v)\in\mathcal{M}\mathcal{R}_p(J;X_1,X_0)$ for each $v\in X_{\gamma,\mu}$. To this end, let $v\in X_{\gamma,\mu}$ be fixed and consider the operator
$$\mathcal{A}(x,D):=\sum_{|\sigma|=4}a_\sigma(\nabla v(x))D^\sigma.$$
For arbitrary but fixed $x_0\in \overline{\Omega}$ we make use of \eqref{eq:appl2a} to calculate the symbol of $\mathcal{A}(x_0,D)$ to the result
$$\mathcal{A}(x_0,\xi)=\sum_{|\sigma|=4}a_\sigma(\nabla v(x_0))\xi^\sigma=
\left(|\xi|^2-\beta_0(\nabla v(x_0)|\xi)^2\right)^2,$$
where $\beta_0:=1/\sqrt{1+|\nabla v(x_0)|^2}$. The Cauchy-Schwarz inequality implies the estimate
$$\mathcal{A}(x_0,\xi)\ge \left(1-\frac{|\nabla v(x_0)|}{\sqrt{1+|\nabla v(x_0)|^2}}\right)^2>0,$$
which is valid for all $\xi\in\mathbb{R}^n$ with $|\xi|=1$. Therefore the symbol of $\mathcal{A}(x_0,D)$ is parameter-elliptic.

Next we verify the Lopatinskii-Shapiro condition (see \cite[Section 2]{DHP2}) for $(\mathcal{A}(x_0,D),\mathcal{B}_j(x_0,D))$, $j\in\{1,2\}$, where $x_0\in \partial\Omega$ is arbitrary but fixed and
$$\mathcal{B}_1(x_0,D)h:=h|_{\partial\Omega}\quad\text{and}\quad \mathcal{B}_2(x_0,D)h:=\nabla h|_{\partial\Omega}\cdot\nu(x_0).$$
To be precise, we have to show that for all $\lambda\in\overline{\mathbb{C}_+}$ and $\xi\in\mathbb{R}^n$ with $(\xi|\nu(x_0))=0$ and $|\xi|+|\lambda|\neq 0$, the only solution $h\in C_0(\mathbb{R}_+)$ of the ODE system
\begin{align}
\begin{split}\label{eq:appl6}
\lambda h(x_n)+\mathcal{A}(x_0,\xi+i\nu(x_0)\partial_n)h(x_n)&=0,\quad x_n>0,\\
\mathcal{B}_j(x_0,\xi+i\nu(x_0)\partial_n)h|_{x_n=0}&=0,\quad j\in\{1,2\},
\end{split}
\end{align}
is $h=0$. It is not difficult to show that
$$\mathcal{A}(x_0,\xi+i\nu(x_0)\partial_n)h=\left(|\xi|^2-\beta_0(\nabla v(x_0)|\xi)^2-\partial_n^2\right)^2h,$$
since $(\nabla v(x_0)|\nu(x_0))=0$ for $v\in X_{\gamma,\mu}$. For the sake of readability, let $b(x_0,\xi):=|\xi|^2-\beta_0(\nabla v(x_0)|\xi)^2$. This yields the representation
$$
\mathcal{A}(x_0,\xi+i\nu(x_0)\partial_n)h=b(x_0,\xi)^2h-2b(x_0,\xi)\partial_n^2 h+\partial_n^4 h.
$$
Furthermore it holds that
$$\mathcal{B}_2(x_0,\xi+i\nu(x_0)\partial_n)h|_{x_n=0}=-\partial_n h|_{x_n=0}.$$
Since $b(x_0,\xi)>0$ for $\xi\neq 0$ it follows that the characteristic polynomial
\begin{equation}\label{eq:appl7}
z^4-2b(x_0,\xi)z^2+b(x_0,\xi)+\lambda=0
\end{equation}
of $\eqref{eq:appl6}_1$ has precisely two roots $z_1,z_2$ with negative real part and two roots $z_3,z_4$ with positive real part (counted by multiplicity). Since we are interested in solutions $h\in C_0(\mathbb{R}_+)$ we may neglect the roots $z_3,z_4$. If $z_1=z_2$, then the general solution of $\eqref{eq:appl6}_1$ is given by
$$h(x_n)=c_1 e^{z_1 x_n}+c_2 x_n e^{z_1 x_n},$$
whereas in case $z_1\neq z_2$, the general solution of $\eqref{eq:appl6}_1$ reads
$$h(x_n)=c_1 e^{z_1 x_n}+c_2 e^{z_2 x_n}.$$
In both cases we may invoke the boundary conditions $\eqref{eq:appl6}_2$ to conclude that $c_1=c_2=0$, hence $h=0$ as desired.

We may apply \cite[Theorem 2.1]{DHP2} to conclude that $A(v)\in\mathcal{M}\mathcal{R}_p(J;X_1,X_0)$ for each $v\in X_{\gamma,\mu}$.
This yields the following result.
\begin{thm}\label{thm:appl2}
Let $n\in\N$, $p,q\in (1,\infty)$ with $4/p+n/q<3$, $\mu_0:=1/p+n/4q+1/4$, $\mu\in (\mu_0,1)$, $\Omega\subset\R^n$ a bounded domain with boundary $\pa\Om\in C^4$. Assume that $h_0\in B_{qp}^{4\mu-4/p}(\Om)$ with $h_0|_{\partial\Omega}=\partial_\nu h_0=0$ at $\partial\Omega$. Then there exists a unique solution $h(t)$ of \eqref{eq:appl2aa} subject to the boundary and initial conditions \eqref{eq:appl2aaa}, with
$$h\in H_{p,\mu}^1(0,T;L_q(\Omega))\cap L_{p,\mu}(0,T;H_q^4(\Omega))\cap C([0,T];B_{qp}^{4\mu-4/p}(\Omega))$$
for each $T\in (0,t^+(h_0))$, where $t^+(h_0)>0$ is the maximal time of existence.

If in addition the solution satisfies
    $$h\in BC\left([\tau,t^+(h_0));B_{qp}^{4\bar{\mu}-4/p}(\Om)\right),$$
for some $\tau\in (0,t^+(h_0))$, $\bar{\mu}\in (\mu,1]$, then $h(t)$ exists globally and the set $\{h(t)\}_{t\ge \tau}$ is relatively compact in $B_{qp}^{4-4/p}(\Om)$. Moreover, the $\om$-limit set
    $$\om(h_0):=\left\{v\in B_{qp}^{4-4/p}(\Om):\ \exists\ t_n\nearrow\infty\ \mbox{s.t.}\ h(t_n,h_0)\to v\ \mbox{in}\ B_{qp}^{4-4/p}(\Om)\right\}$$
is nonempty, connected and compact.
\end{thm}

\subsection{Willmore flow - The graph case}

We consider a family of hypersurfaces  $\{\Gamma(t)\}_{t\ge 0}\subset\mathbb{R}^{n+1}$ in the same setting as in Subsection \ref{subsec:surfdiff} but this time we assume that the evolution of $\Gamma(t)$ is induced by the Willmore flow. To be precise we have
\begin{equation}\label{eq:appl8}
V_{\Gamma(t)}=-\Delta_{\Gamma(t)}H_{\Gamma(t)}+H_{\Gamma(t)}
\left[\frac{1}{2}H_{\Gamma(t)}^2- \operatorname{tr}L_{\Gamma(t)}^2\right],
\end{equation}
where $V_{\Gamma(t)}$, $\Delta_{\Gamma(t)}$, $H_{\Gamma(t)}$ have the same meaning as in Subsection \ref{subsec:surfdiff} and $L_{\Gamma(t)}$ denotes the \emph{Weingarten tensor}; see \cite[Subsection 2.2]{PrSi13} for its definition. It follows from \cite[Subsections 2.3 \& 2.8]{PrSi13} that in the graph case the formula
$$\operatorname{tr}L_{\Gamma}^2=-(\delta^{ij}-\beta^2\partial_i h\partial_j h)(\partial_i\partial_j\nu_\Gamma|\nu_{\Gamma})$$
holds true, where we made use of sum convention. Note that the outer unit normal field $\nu_{\Gamma(t)}$ on $\Gamma(t)$ is given by $\nu_\Gamma=\beta(-\nabla h,1)^{\sf T}$ and $\beta=1/\sqrt{1+|\nabla h|^2}$.

From the representations in Subsection \ref{subsec:surfdiff} we obtain
$$H_\Gamma^3=(\delta^{ij}-\beta^2\partial_ih\partial_jh)^3\beta^3(\partial_i\partial_j h)^3,$$
and
$$H_\Gamma\operatorname{tr}L_\Gamma^2=-(\delta^{ij}-\beta^2\partial_ih\partial_jh)^2\beta
\partial_i\partial_j h(\partial_i\partial_j\nu_\Gamma|\nu_{\Gamma}).$$
Concerning the term $\partial_i\partial_j\nu_\Gamma$ we compute
$$\partial_i\partial_j\nu_\Gamma=(\partial_i\partial_j\beta)(-\nabla h,1)^{\sf T}-\partial_j\beta(\partial_i\nabla h,0)^{\sf T}-\partial_i\beta(\partial_j\nabla h,0)^{\sf T}-\beta(\partial_i\partial_j\nabla h,0)^{\sf T},$$
for $i,j\in\{1,\ldots,n\}$. Furthermore, since $\partial_j\beta=-\beta^3(\partial_j\nabla h|\nabla h)$, we obtain
$$\partial_i\partial_j\beta=-3\beta^2(\partial_i\beta)(\partial_j\nabla h|\nabla h)
-\beta^3(\partial_i\partial_j\nabla h|\nabla h)-\beta^3(\partial_j\nabla h|\partial_i\nabla h),$$
for $i,j\in\{1,\ldots,n\}$. The above computations show that
\begin{multline*}
H_\Gamma\left[\frac{1}{2}H_{\Gamma}^2- \operatorname{tr}L_{\Gamma}^2\right]\\
=\sum_{|\sigma|=3,|\tau|=2}\tilde{b}_{\sigma\tau}(\nabla h)D^\sigma h D^\tau h
+\sum_{|\sigma|=|\tau|=|\chi|=2}\tilde{c}_{\sigma\tau\chi}(\nabla h)D^\sigma h D^\tau h D^\chi h,
\end{multline*}
where the coefficients $\tilde{b}_{\sigma\tau}$ and $\tilde{c}_{\sigma\tau\chi}$ are smooth in $\nabla h$. Therefore the height function $h$ satisfies the partial differential equation \eqref{eq:appl2aa} with the same leading coefficient $a_\sigma$ and some modified and smooth coefficients $b_{\sigma\tau}$, $c_{\sigma\tau\chi}$. Imposing the boundary conditions $h|_{\partial\Omega}=\partial_\nu h=0$ at $\partial\Omega$ it follows that Theorem \ref{thm:appl2} also holds for the Willmore flow.

%{\diver}_{x'}(\beta\nabla _{x'} h)=\beta[\Delta_{x'}h-\beta^2(\nabla^2_{x'} h\nabla_{x'} h|\nabla_{x'} h)]
%
%\Delta_{x'}\varphi-\beta^2(\nabla_{x'}^2\varphi\nabla_{x'} h|\nabla_{x'} h)
%-\beta^2[\Delta_{x'} h-\beta^2(\nabla_{x'}^2 h\nabla_{x'} h|\nabla_{x'} h)](\nabla_{x'} h|\nabla_{x'}\varphi)
%%%%%%%%%%%%%%%%%%%%%%%%%%%%%%%%%
\bibliographystyle{amsplain}

\end{document}